# On Nilpotence of a Kind of Circulant Matrices over $Z_p$

Wei Wang

College of Information and Engineering, Tarim University

Alar, Xinjiang 843300

**Abstract.** We investigate the nilpotence of a kind of circulant matrices $T_{n,m}$ over field $Z_p$, where $T_{n,m} = \sum_{i=0}^{m-1} S_n^i$ and $S_n$ is the fundamental circulant matrix of order $n$. The necessary and sufficient condition on $n$ and $m$ for determining nilpotence of $T_{n,m}$ over $Z_p$ is presented. Moreover, we obtain a formula for nilpotent index of $T_{n,m}$ when the condition is satisfied. As an application, we give a complete solution of a conjecture of C.Y.Zhang.

## 1. Introduction

In this paper, motivated by a conjecture of C.Y.Zhang(see[1], and[2]for a partial solution), we investigate the nilpotence of a special kind of circulant matrices $T_{n,m}$ over $Z_p$, where $Z_p$ denotes residue class field modulo $p$. For $n, m > 0$, define $T_{n,m}$ as follows.

$$T_{n,m} = \sum_{i=0}^{m-1} S_n^i = I_n + S_n + S_n^2 + \cdots + S_n^{m-1},$$

where $I_n, S_n$ denote, respectively, identity matrix and the fundamental circulant matrix of order $n$.

Any integer matrix can be regarded as a matrix over field $Z_p$, and two integer matrices $A = (a_{ij})$, $B = (b_{ij})$ of the same shape are equal over $Z_p$ whenever $A \equiv B (\mod p)$, i.e., $a_{ij} \equiv b_{ij} (\mod p)$ for all $i, j$. See[3] for a general view of matrix operations over $Z_p$.

A integer matrix $M$ is called *nilpotent* over $Z_p$ if $M^k \equiv 0 (\mod p)$ for some $k > 0$. For a nilpotent matrix $M$ over $Z_p$, the minimal $k$ with $M^k \equiv 0 (\mod p)$ is called the *nilpotent index* of $M$. Similar concepts can be defined over ring $Z_q$ when $q$ is not necessarily prime. We remark that $M$ is nilpotent over $Z_q$ if and only if $M$ is nilpotent over $Z_p$ whenever $p$ is a prime factor of $q$. For notations and terminology not defined here, we refer to [4,5].

The main result of this paper is stated in the following theorem.

**Theorem 1:** Let $n, m > 0$ and $p$ be a prime. Then $T_{n,m}$ is nilpotent over $Z_p$ if and only if (i) $p \mid m$ and (ii) $n \mid mp^k$ for some $k \geq 0$. Moreover, the nilpotent index of $T_{n,m}$ over $Z_p$ is given by $\left\lceil \dfrac{p^a}{p^b - 1} \right\rceil$, where $a, b$ are the largest integers with $p^a \mid n$, $p^b \mid m$, respectively.

Notice $b > 0$ by condition (i). Proof for Theorem 1 is postponed in Sec.3 after we establish some preliminary results in Sec.2. Based on Theorem 1, it is not difficult to get the following corollary, which essentially gives a complete solution of a conjecture of C.Y.Zhang.

**Corollary 1:** Let $n, m$ be positive integers with $m > 1$. Then $T_{n,m}$ is nilpotent over $Z_m$ if and only if either (i) Both $m$ and $n$ are some powers of the same prime, or (ii) $m$ has at least two different prime factors and $n \mid m$. Moreover, the nilpotent index of $T_{n,m}$ over $Z_m$ is $\leq n$ in either cases.

## 2. Preliminaries

**Lemma 1:** Let $m = dm^*, n = d^q n^*$ with $q > 0$, $(d, m^*) = (d, n^*) = 1$ and $n^* \mid m^*$. For any integer $c$, the following congruence equation with $q$ variables has exactly $\dfrac{(m^*)^q}{n^*}$ integer solutions with $0 \leq x_i < m$ for $0 \leq i < q$.

$$x_0 + dx_1 + \cdots + d^{q-1} x_{q-1} \equiv c \pmod{n}. \tag{1}$$

**Proof:** By induction on $q$. When $q = 1$, (1) becomes

$$x_0 \equiv c \pmod{n}. \tag{2}$$

Notice $n \mid m$. Hence (2) has $\dfrac{m}{n}$ incongruent solutions modulo $m$.

$$x_0 \equiv c + jn \pmod{m}, 0 \leq j < \dfrac{m}{n}. \tag{3}$$

Therefore, the lemma is true when $q = 1$.

Next we assume that the lemma is true for $q = k$, and we proceed to check it for $q = k+1$. We need to count the number of integer solutions with $0 \leq x_i < m$ of the following equation.

$$x_0 + dx_1 + \cdots + d^k x_k \equiv c \pmod{n}. \tag{4}$$

Since $d \mid n$, any solution of (4) must satisfy $x_0 \equiv c \pmod{d}$, which has exactly $\dfrac{m}{d} = m^*$ solutions with $0 \leq x_0 < m$. Let $x_0^*$ be such a solution. By substituting $x_0 = x_0^*$ to (4), we have

$$dx_1 + d^2 x_2 + \cdots + d^k x_k \equiv c - x_0^* \pmod{n}. \tag{5}$$

Let $c' = \dfrac{c - x_0^*}{d}$ and $n' = \dfrac{n}{d} = d^k n^*$. Then (5) is equivalent to

$$x_1 + dx_2 + \cdots + d^{k-1} x_k \equiv c' \pmod{n'}. \tag{6}$$

By the induction hypothesis, (6) has exactly $\dfrac{(m^*)^k}{n^*}$ solutions with $0 \leq x_i < m$ for $1 \leq i < k+1$. Thus, (4) has exactly $m^* \dfrac{(m^*)^k}{n^*} = \dfrac{(m^*)^{k+1}}{n^*}$ solutions with $0 \leq x_i < m$ for $0 \leq i < k+1$, which completes the induction proof.

**Lemma 2:** Let $k \geq 0$ and $R$ be a commutative ring with characteristic $p$. Then the map $\varphi: R \to R$ defined by $\varphi(x) = x^{p^k}$ is a ring homomorphism.

We remark that all circulant matrices of order $n$ over $Z_p$ consist a commutative ring with characteristic $p$.

## 3. Proof of Theorem 1

For simplicity, we will omit subscripts of $T_{n,m}$ and $S_n$.

**(Necessity)** Let $\xi = (1,1,\cdots,1)^T$ be an $n$-dimensional vector. The nilpotence of $T$ over $Z_p$ implies
$$T^k \xi = m^k \xi \equiv 0 (\bmod p), \text{ for some } k > 0.$$
Thus $m^k \equiv 0 (\bmod p)$ and hence $m \equiv 0 (\bmod p)$, i.e., $p \mid m$.

Next, by the identity $(I + S + \cdots + S^{m-1})(I - S) = I - S^m$, the nilpotence of $T$ implies $I - S^m$ is nilpotent, i.e., $(I - S^m)^{p^k} \equiv 0 (\bmod p)$, for some $k \geq 0$. Hence, by Lemma 2, we have
$$I - S^{mp^k} \equiv 0 (\bmod p), \text{ for some } k \geq 0,$$
which holds only if $I - S^{mp^k} = 0$ since any power of $S$ is a permutation matrix. Therefore, $n \mid mp^k$.

**(Sufficiency)** Let $n = p^a n^*$, $m = p^b m^*$. By division with remainder, $a = bq + r$ with $0 \leq r < b$. We will prove this part by directly showing (A) $T^{\left\lceil \frac{p^a}{p^b-1} \right\rceil - 1} \not\equiv 0 (\bmod p)$ and (B) $T^{\left\lceil \frac{p^a}{p^b-1} \right\rceil} \equiv 0 (\bmod p)$. We may assume $q > 0$ since otherwise one can easily show $T \equiv 0 (\bmod p)$.

It is simple to check
$$\left\lceil \frac{p^a}{p^b - 1} \right\rceil = \left\lceil \frac{p^r(p^{bq}-1)+p^r}{p^b-1} \right\rceil = p^r(1 + p^b + \cdots + p^{b(q-1)}) + 1.$$

To show (A), we do the following calculation.
$$\left(\sum_{i=0}^{m-1} S^i\right)^{p^r(1+p^b+\cdots+p^{b(q-1)})} = \left(\left(\sum_{i=0}^{m-1} S^i\right)^{1+p^b+\cdots+p^{b(q-1)}}\right)^{p^r} \equiv \left(\prod_{j=0}^{q-1}\sum_{i=0}^{m-1} S^{ip^{bj}}\right)^{p^r} \quad \text{(by Lemma 2)}$$

$$= \left(\sum_{0 \leq i_0, i_1, \cdots, i_{q-1} < m} S^{i_0 + i_1 p^b + \cdots + i_{q-1} p^{b(q-1)}}\right)^{p^r}$$

$$\equiv \sum_{0 \leq i_0, i_1, \cdots, i_{q-1} < m} S^{p^r(i_0 + i_1 p^b + \cdots + i_{q-1} p^{b(q-1)})} \quad \text{(by Lemma 2)}$$

$$= \frac{(m^*)^q}{n^*} \sum_{c=0, p^r \mid c}^{n-1} S^c \quad \text{( by Lemma 1, } S^n = I \text{ and the fact explained bellow)}$$

$\not\equiv 0 \pmod{p}$  (since $(m^*, p) = (n^*, p) = 1$)

In the last equality we used the fact that $p^r(x_0 + x_1 p^b + \cdots + x_{q-1} p^{b(q-1)}) \equiv c \pmod{p^{bq+r} n^*}$ is equivalent to $x_0 + x_1 p^b + \cdots + x_{q-1} p^{b(q-1)} \equiv \dfrac{c}{p^r} \pmod{p^{bq} n^*}$ if $p^r \mid c$ and has no solutions otherwise.

We are left to check (B). By division with remainder, each $i$ with $0 \leq i < m = p^b m^*$ can be written uniquely in the form $i = kp^r + j$ with $0 \leq k < p^{b-r} m^*$ and $0 \leq j < p^r$. And vice versa. Hence,

$$\sum_{c=0, p^r \mid c}^{n-1} S^c \sum_{i=0}^{m-1} S^i = \sum_{c=0, p^r \mid c}^{n-1} S^c \sum_{k=0}^{p^{b-r} m^* - 1} S^{kp^r} \sum_{j=0}^{p^r - 1} S^j$$

$$= p^{b-r} m^* \sum_{c=0, p^r \mid c}^{n-1} S^c \sum_{j=0}^{p^r - 1} S^j \equiv 0 \pmod{p}$$

The second equality follows from the fact $\sum_{c=0, p^r \mid c}^{n-1} S^c S^{p^r} = \sum_{c=0, p^r \mid c}^{n-1} S^c$.